\documentclass[3p]{elsarticle}
\usepackage{amsfonts}
\usepackage{amsmath}
\usepackage{theorem}
\usepackage{amssymb}
\usepackage[T1]{fontenc}
\usepackage{lmodern}
\usepackage{graphicx}
\usepackage{hyperref}
\usepackage[noend]{algpseudocode}
\usepackage{xspace}
\usepackage{algorithm}
\usepackage{sectsty}



\algrenewcommand\algorithmicrequire{\textbf{Input:}}
\algrenewcommand\algorithmicensure{\textbf{Output:}}

\theoremstyle{break}
\newtheorem{theorem}{Theorem}[section]

\theorembodyfont{\rm}

\newtheorem{experiment}{Experiment}[section]

\newcommand\spann{\mathop\mathrm{span}\nolimits}

\newcommand\diag{\mathop\mathrm{diag}\nolimits}

\newcommand{\cnn}{\mathbb{C}^{n\times n}}

\newcommand{\calx}{\mathcal{X}}

\newcommand{\ilambda}{I_{\lambda}}

\newcommand{\imag}{i}

\newcommand{\dz}{\mathrm{d}z}
\newcommand{\feast}{FEAST\xspace}
\newcommand{\RR}{Rayleigh--Ritz\xspace}
\newcommand{\RRm}{\RR method\xspace}

\newcommand{\ie}{i.\,e.,\;}
\newcommand{\eg}{e.\,g.,\;}

\newcommand{\lmin}{\underline{\lambda}}
\newcommand{\lmax}{\overline{\lambda}}
\newcommand{\abs}[1]{\left\lvert #1 \right\rvert}
\newcommand{\norm}[1]{\left\| #1 \right\|}
\newcommand{\figwidth}{79mm}

\title{Dissecting the FEAST algorithm for generalized eigenproblems}
\author[buw]{Lukas Kr{\"a}mer\fnref{fn1}}
\ead{lkraemer@math.uni-wuppertal.de}
\author[juel]{Edoardo Di~Napoli\corref{cor}\fnref{fn3}}
\ead{dinapoli@aices.rwth-aachen.de}
\author[buw]{Martin Galgon\fnref{fn1}}
\ead{galgon@math.uni-wuppertal.de}
\author[buw]{Bruno Lang\fnref{fn1}}
\ead{lang@math.uni-wuppertal.de}
\author[rwth]{Paolo Bientinesi\fnref{fn2}}
\ead{pauldj@aices.rwth-aachen.de}

 \fntext[fn1]{These authors gratefully acknowledge support by the Bundesministerium
  f\"ur Bildung und Forschung within the project ``ELPA --
  Hochskalierbare Eigenwert-L\"oser f\"ur
  Petaflop-Gro{\ss}anwendungen'', F\"orderkennzeichen 01IH08007B.}
\fntext[fn2]{Financial support from the Deutsche
  Forschungsgemeinschaft (German Research Association) through grant
  GSC 111 is gratefully acknowledged.}
\fntext[fn3]{Financial support by the Volkswagen Foundation through the fellowship ``Computational Sciences'' is gratefully acknowledged.}

\cortext[cor]{Corresponding author.}
\address[buw]{Bergische Universit\"at Wuppertal,
     Fachbereich C -- Mathematik und Naturwissenschaften,
     Gau\ss str. 20,
     42119 Wuppertal, Germany
     }

\address[rwth]{
     RWTH Aachen,
     AICES,
     Schinkelstr. 2,
     52062 Aachen, 
     Germany
     }

\address[juel]{
  Forschungszentrum J\"ulich,
  Institute for Advanced Simulation, 
  J\"ulich Supercomputing Centre, 
  52425 J\"ulich, 
  Germany
}

\begin{document}

\begin{abstract}
We analyze the \feast method
for computing selected eigenvalues and eigenvectors of large sparse
matrix pencils. After establishing the close connection between \feast
and the well-known \RRm, we identify several critical issues that
influence convergence and accuracy of the solver:
the choice of the starting vector space, the stopping criterion, 
how the inner linear systems impact the quality of the solution,  
and the use of \feast for computing eigenpairs from multiple intervals.
We complement the study with numerical examples, and hint at possible 
improvements to overcome the existing problems.
\end{abstract}


\maketitle

\section{Introduction}\label{sec:Intro}

In 2009, Polizzi introduced the \feast\ solver for generalized
Hermitian definite eigenproblems~\cite{polizzi:2009}.
\feast\ was conceived as an algorithm for electronic structure
calculations, and then evolved into a general purpose solver. 
In this paper we describe the mathematical structure of the algorithm
and conduct an investigation of its robustness and accuracy.

\feast\ belongs to a family of iterative solvers based on the contour
integration of a density-matrix representation of quantum mechanics; the
result of the integration is a subspace projector that plays a central
role in a \RRm. The prominent member in this family was developed by
Sakurai and Sugiura in 2003~\cite{sakurai_sugiuara:2003}; 
this work, which inspired a number of generalizations and 
high-performance implementations~\cite{sakurai2, sakurai3}, 
targets non-Hermitian eigenproblems. By contrast, \feast\ promises to
deliver performance on sparse Hermitian problems.  Such problems can
be solved by a number of alternative packages like ARPACK~\cite{arug:1998}
and TRLan~\cite{trlan} and the solver implemented in PARSEC~\cite{parsec1, parsec2}.

Since in ab-initio electronic calculations typically one is interested
in the lowest part of the eigenspectrum, we investigate \feast's
behavior for the computation of a subset of eigenpairs lying inside a
given interval. Additionally, we study its strengths and weaknesses
when a large portion or all of the spectrum is sought
after. In our analysis, we use a number of matrices from practical
applications.  From our experiments, we found that while for specific
scenarios \feast\ is accurate and reliable, in general it lacks
robustness.

Our analysis touches upon three main features of the solver: 1)
critical input parameters, 2) the stopping criterion, and 3) the
quality of the results.

\begin{itemize}
\item[1)] 
In addition to a
search interval, \feast's interface requires the user to specify the
number of eigenvalues present within the interval.  Although in some
ab-initio simulations this number can be accurately estimated, in
general it is not possible to obtain it cheaply; since the completeness
of the computed eigenpairs greatly depends on it, this initial guess
is critical.  The user can also specify a starting vector base which \feast uses to initialize
 the solver: while by default \feast uses a
random set of vectors, the convergence rate of the solver greatly
depends on the actual choice.  We elaborate on how different starting
bases affect convergence speed and robustness.
\item[2)] 
The original implementation of \feast employs a stopping criterion
based on monitoring relative changes in the sum of all computed Ritz
values. We identify cases where this criterion does not reflect the
actual convergence and propose an alternative criterion based on
per-eigenpair residuals.
\item[3)] 
  We analyze the quality of the solution computed by \feast by
  means of residuals and orthogonality. On the one hand we found that 
  the achievable residuals are affected by the accuracy of the linear solver 
  used within the algorithm. On the other hand, 
  we distinguish between local and global orthogonality, 
  \ie orthogonality among eigenvectors that were computed with a single interval,
  and among vectors from separate intervals, respectively.
  While the local orthogonality is guaranteed by the \RRm, 
  depending on the spectrum of the eigenproblem 
  the global orthogonality might suffer.
\end{itemize}

The paper is organized as follows.  In Section~\ref{sec:FEAST-RR} we
illustrate the underlying mathematical structure of the algorithm.
Section~\ref{sec:study} contains experiments and analysis concerning
the different aspects of the solver.
Section~\ref{feast_in_context:sec} examines the suitability of \feast
as a building-block for a general purpose solver, working on multiple
intervals to compute a larger portion of the spectrum. We conclude in
Section~\ref{sec:concl} suggesting improvements to broaden
\feast's applicability.

\section{FEAST and the \RRm}\label{sec:FEAST-RR}

Let us consider the generalized eigenproblem $A x = \lambda B x$, 
with Hermitian matrix $A \in \cnn$ and Hermitian positive definite $B\in\cnn$.
The objective is to compute the eigenpairs whose eigenvalues lie
in a given interval $\ilambda = [ \lmin, \lmax ]$.
Since \feast\ is an instantiation of the \RRm, we start with a short review.

In the following, in order to denote an eigenpair $(\lambda,x)$ 
with $\lambda\in \ilambda$, we will sloppily say that the eigenpair is 
in the interval $\ilambda$.

\subsection{\RR\ theorem}

We present the \RRm\ in its orthonormal version, which ensures a minimal
residual. The method relies on the following theorem.

\begin{theorem}[\RR, \cite{stewart:2001}]\label{RR:thm}
  Let $\mathcal{U}$ be a subspace containing an eigenspace
  $\mathcal{X} \subset \mathcal{U}$ of the generalized eigenproblem 
  $A x = \lambda B x$.
  Let $U$ be a basis of vectors for
  $\mathcal{U} = \spann(U)$, $U^I$ a left inverse of $U$, and $A_U = U^I A
  U$, $B_U = U^I B U$ the so-called Rayleigh quotients for $A$ and
  $B$. If $(\Lambda,W)$ 
  (${\Lambda}=\diag(\lambda_1,\lambda_2,\ldots)$)
  are \emph{primitive Ritz pairs} of the
  reduced problem, \ie $A_U W = B_U W \Lambda$, then
  $(\Lambda,X)$ are {\em Ritz pairs} for the original eigenproblem
  with $X = U W$ and $\spann(X) = \mathcal{X}$.
\end{theorem}
In a neighborhood of the exact solutions for the primitive Ritz pair,
we expect that the \RR\ theorem also applies approximately,
leading to the following strategy.
\begin{enumerate}
  \item Find a suitable basis $U$ for $\mathcal{U}$.
  \item Compute the Rayleigh quotients $A_U = U^I \! A U$,
        $B_U = U^I \! B U$.
  \item Compute the primitive Ritz pairs $(\widetilde{\Lambda},\widetilde{W})$
        of $A_U W = B_U W \Lambda$.
  \item Return the approximate Ritz pairs $(\widetilde{\Lambda},U\widetilde{W})$
        of $A X = B X \Lambda$.
  \item Check convergence criterion; if not satisfied, go back to
        Step~1.
\end{enumerate}

Let us point out that obtaining an accurate approximation of
$(\Lambda,X)$ is not an obvious consequence of computing primitive
Ritz pairs.  One must ensure that both Ritz values and the
corresponding Ritz vectors converge to the desired eigenpairs.
Conversely, it must hold that each eigenpair in $\ilambda$ corresponds
to a primitive Ritz pair in the same interval;
see~\cite{stewart:2001}.

\subsection{The algorithm}

The \feast algorithm implements the above \RRm, with a particular choice 
for computing $U$:
\begin{equation}\label{eq:u-integral}
  U := \frac{1}{2\pi\imag} \int_{\mathcal{C}}\dz (zB-A)^{-1} B Y,
\end{equation}
where $\mathcal{C}$ is a curve in the complex plane enclosing the
selected interval $\ilambda$.  The expression $(zB-A)^{-1} B$ is
normally referred to as the eigenproblem's {\em resolvent};
Formula~\eqref{eq:u-integral} can be interpreted as the projection of
the set of vectors $Y$ onto a subspace $\mathcal{U}$ containing the
eigenspace.  Pseudo-code for \feast is provided in
Algorithm~\ref{basicfeast:alg}.
\begin{algorithm}[htbp]
\caption{Skeleton of the \feast algorithm}
\label{basicfeast:alg}
\begin{algorithmic}[1]
  \Require An interval $\ilambda = \left[\lmin,\lmax\right]$ and
           an estimate $\widetilde{M}$ of the number of eigenvalues in $\ilambda$.
  \Ensure $\hat{M} \le \widetilde{M}$ eigenpairs in $\ilambda$.
  \State\label{integral:line}%
         Choose $Y\in \mathbb{C}^{n\times \widetilde{M}}$ of rank $\widetilde{M}$ 
         and compute
         $U := \frac{1}{2\pi\imag} \int_{\mathcal{C}}\dz\ (zB-A)^{-1} B\, Y$;
  \State \label{RR_line1:line}%
          Form the Rayleigh quotients $A_U:= U^H\!AU, \; B_U:=U^H\!BU$;
  \State \label{RR_line2:line}%
          Solve the size-$\widetilde{M}$ generalized eigenproblem
          $A_U\widetilde{W} = B_U\widetilde{W}\widetilde{\Lambda}$;
  \State Compute the approximate Ritz pairs
         $( \widetilde{\Lambda},  \widetilde{X}:=U\cdot\widetilde{W} )$;
  \State If convergence is not reached\label{stopping:line}
         then go to Step~\ref{integral:line}, with $Y := \widetilde{X}$.
\end{algorithmic}
\end{algorithm}

Having established the close connection between \feast and the \RRm,
the next section is devoted to
the theoretical framework for the resolvent and 
the contour integration~(\ref{eq:u-integral}).

\subsection{Integrating the resolvent}%
   \label{sec:IntegratingTheResolvent}

   In this section we define the concept of resolvent and illustrate
   its functionality within the \RRm. The objective is to show that
   the subspace $U$ is approximated by the integral operator in 
   Equation~\eqref{eq:u-integral}.
   We recall that
   a generalized Hermitian definite eigenproblem 
   has $n$ real eigenvalues $\lambda_{1}$, \ldots, $\lambda_{n}$ and $B$-orthonormal
   eigenvectors $x_{1}$, \ldots, $x_{n}$.

Let us consider a single eigenpair $(\lambda_{k}, x_{k})$, and
let $z \in \mathbb{C}$. For $z \not= \lambda_{i}, i \in \{1, \dots, n\}$,
\[
   B^{-1} (z B - A) x_{k} = (z - \lambda_{k}) x_{k},
\]
and
\begin{equation}
    (z B - A)^{-1} B x_{k}
  = (B^{-1}(z B - A))^{-1} x_{k}
  = (z - \lambda_{k})^{-1} x_{k}
  .
  \label{eq:pole}
\end{equation}
Define the resolvent operator $G(z)$ as
\[
  G( z ) := (z B - A)^{-1} B  
\]
and let $\mathcal{C}_{k}$ be a closed curve in the complex plane enclosing 
only the eigenvalue $\lambda_{k}$.
Thus, the integral
\[
  \frac{1}{2 \pi \imag} \int_{\mathcal{C}_{k}} \dz \, G(z) x_{k}
\]
equals the residue of $G(z) x_{k}$ localized at the pole in $\lambda_{k}$.
Using Equation~\eqref{eq:pole}, we obtain
\begin{equation}
    \frac{1}{2 \pi \imag}
       \int_{\mathcal{C}_{k}} \dz \, (z B - A)^{-1} B x_{k}
  = \frac{1}{2 \pi \imag}
       \int_{\mathcal{C}_{k}} \frac{\dz}{z - \lambda_{k}} x_{k}
  = \frac{1}{2 \pi \imag} 2 \pi \imag x_{k}
  = x_{k}
  .
  \label{eq:int-k}
\end{equation}
By contrast, the residue around any other pole returns 0.

Now let $\mathcal{C}$ be a curve enclosing a subset
$\{\lambda_{k} : k \in I \}$ of the eigenvalues.
Combining Equations~\eqref{eq:int-k} for $k\in I$ and
splitting the path integral over $\mathcal{C}$ into a sum of integrals
over closed curves $\mathcal{C}_{k}$ containing just one eigenvalue
$\lambda_{k}$ each, one obtains
\begin{equation}
  \label{eq:int2}
    \frac{1}{2 \pi \imag}
      \int_{\mathcal{C}} \dz \,G(z) x_j
  = \sum_{k \in I}
      \frac{1}{2 \pi \imag} \int_{\mathcal{C}_k} \dz \,G(z) x_j
  = \sum_{k \in I}
      \delta_{k,j} x_{j}
  = \left\{
	\begin{array}{ll}
	   x_{j}, & \mbox{if $j \in I$} \\
         0,     & \mbox{otherwise}
      \end{array}
    \right.
  .
\end{equation}

So far we have shown how the projection operator acts on 
the full space of eigenvectors. This is a
well known property of the resolvent of an eigenproblem~\cite[Ch.~3]{saad:2003}. 
This property can also be visualized by combining the
$B$-orthogonal eigenvectors $x_{k}$, $k \in I$, 
into an $n$-by-$| I |$ matrix $X$, and by comparing~\eqref{eq:int2} with
the application of the projector 
\[
  Q = XX^H\ B = \sum_{k \in I} x_k x_k^H\ B \quad {\rm with} \quad Q^2 = Q 
\] 
to an eigenvector $x_{j}$,
\[ 
    Q x_{j}
  = \sum_{k \in I} (x_k x_k^H) B x_j
  = \sum_{k \in I} x_k \delta_{k,j}
  = \left\{
	\begin{array}{ll}
	   x_{j}, & \mbox{if $j \in I$} \\
         0,     & \mbox{otherwise}
      \end{array}
    \right.
  .
\] 
One then concludes that the operators 
$  \frac{1}{2 \pi \imag} \int_{\mathcal{C}} \dz \, G(z)$ 
and $Q$, when applied to a set of eigenvectors, produce the same results.

Let $Y = \{y_1,y_2,\dots, y_M\}$, then
\[
    \frac{1}{2\pi\imag} \int_{\mathcal{C}}\dz\ (zB-A)^{-1}\ B Y
  = Q Y
  = X X^H B Y
\]
projects each $y_j$ onto the eigenspace $\mathcal{X} = \spann(X)$.\footnote{
In Polizzi's original paper, $U$ equals $XX^HY$ instead of
$XX^HBY$. This is correct only when $Y$ is chosen to be a random
set of vectors, since $BY$ does not alter the random nature of $Y$.
Thus, this equality is only valid in the first iteration of the solver.}
In this sense, the matrix $U$ computed in Algorithm~\ref{basicfeast:alg}
is a reasonable attempt to fulfill the requirements of the \RR theorem.
In practice, the integral~\eqref{eq:u-integral} must be evaluated numerically, using
a scheme such as Gau\ss--Legendre; for details we refer to
the original publication~\cite{polizzi:2009} and to the literature on
numerical integration, \eg~\cite{davis_rabinowitz:1984}.

A numerical integration scheme leads to an approximation
\[
   \hat U
   \approx
   \frac{1}{2\pi\imag}
   \sum_{k=1}^{m} w_{k} (z_{k} B-A)^{-1}\ B Y
\]
where the points $z_{k}$ lie on the curve $\mathcal{C}$.
The accuracy of the approximation, as well as the computational complexity,
are determined by the number of integration points.
In practice, Gaussian rules with 7--10 nodes already achieve
satisfactory results. 
For each integration point, a linear system $( z_{k} B - A ) \hat U_k = BY$
of size $N$ with $\widetilde{M}$ right-hand sides must be solved.

If the curve $\mathcal{C}$ is chosen to be symmetric with respect to the 
real axis (e.g., a circle or ellipse), then considerable 
computational savings are possible.
In this case the numerical integration must cover only the half-curve in
the upper complex half-plane due to the symmetry
$G(\bar{z}) = G^{H}(z)$ \cite{polizzi:2009}.

\section{Analysis and experiments for the \feast algorithm}\label{sec:study}

In this section we discuss the issues arising when employing \feast to compute the eigenpairs in a single interval: size and choice of the search space, stopping criteria and impact of the linear solver on the accuracy of the solution.

\subsection{Size of the search space $\mathcal{U}$}\label{searchspace:sec}

Section~\ref{sec:IntegratingTheResolvent} shows that the size of the 
eigenspace $\mathcal{X}$ is determined by the number $M$ of columns of $U$ 
which in turn corresponds to the number of columns of $Y$.
This number is equivalent to the number $M$ of eigenvalues (counting
multiplicities) lying in the given interval $\ilambda$.
In practice, $M$ is not known a priori, and an estimate
$\widetilde{M}$ has to be used instead.
In the following, we discuss the consequences of the
cases $\widetilde{M} > M$ and $\widetilde{M} < M$.

\textbf{Case $\widetilde{M}>M$.}
If the estimate is larger than the actual number of eigenvalues
in $\ilambda$, then the spectral projector $Q$ has rank $M$, and thus
the $n$-by-$\widetilde{M}$ matrix $U=QY$ is rank-deficient and does not
have a full rank left inverse. 
As a direct consequence, $B_U$ is rank deficient, and  
the eigenpairs of the reduced eigenproblem might be incomplete and 
not necessarily $B$-orthogonal.

In the original implementation of \feast, the rank deficiency is
detected by measuring the ``positive definiteness'' of $B_U$ through a
Cholesky decomposition.  A drawback of
this approach is that it does not indicate which of the columns of $U$
are linearly independent, forcing one to select a subset 
arbitrarily. This, in turn, might result in the presence of
spurious eigenpairs.  A reliable, but more expensive, approach
consists of computing an SVD or a (rank-revealing) QR decomposition of
the matrix $U$ \cite{golubvanloan:1996}.
Since the rank of $U$ equals the number of eigenvalues in $\ilambda$,
such a rank-revealing decomposition can be used to safely restart
the process with $Y = \bar U X$, where $\bar U$ includes the linearly independent columns of $U$.

\textbf{Case $\widetilde{M}<M$.}
The space spanned by $U$ does not contain the whole
eigenspace corresponding to the eigenvalues within the integration
contour, and therefore the vectors generated by integrating the resolvent fail to span $\mathcal{X}$.

The following experiment illustrates the behavior of \feast for different numbers $\widetilde{M}$.

\begin{figure}[htb]
   \begin{center}
      \includegraphics[width=\figwidth]{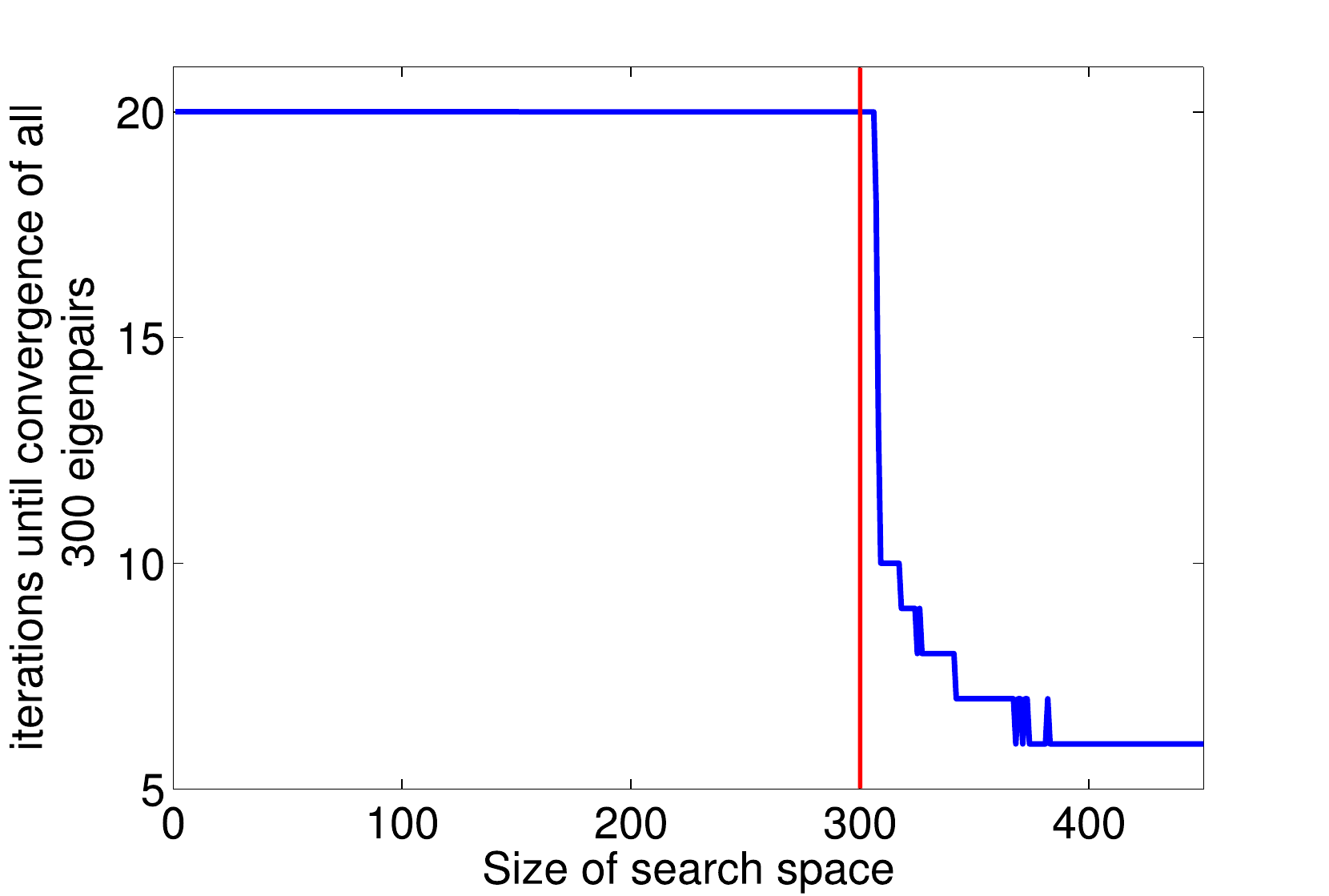}
      \quad
      \includegraphics[width=\figwidth]{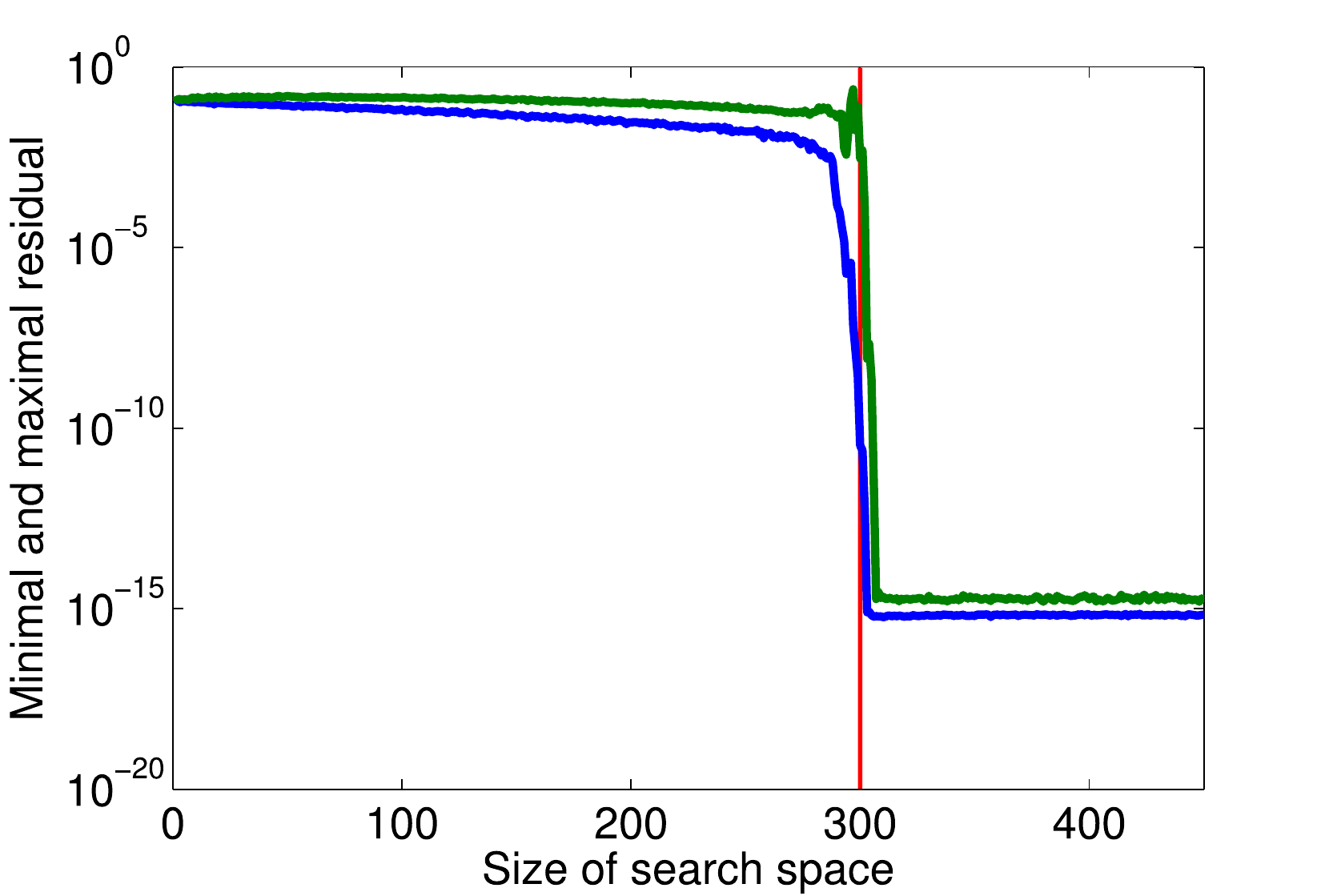}
   \end{center}
   \caption{Left: Number of necessary iterations.
            Right: Minimal (lower line) and maximal (upper line) residual.}
   \label{no_iters-and-minmaxres:fig}
\end{figure}

\begin{figure}[htb]
   \begin{center}
      \includegraphics[width=\figwidth]{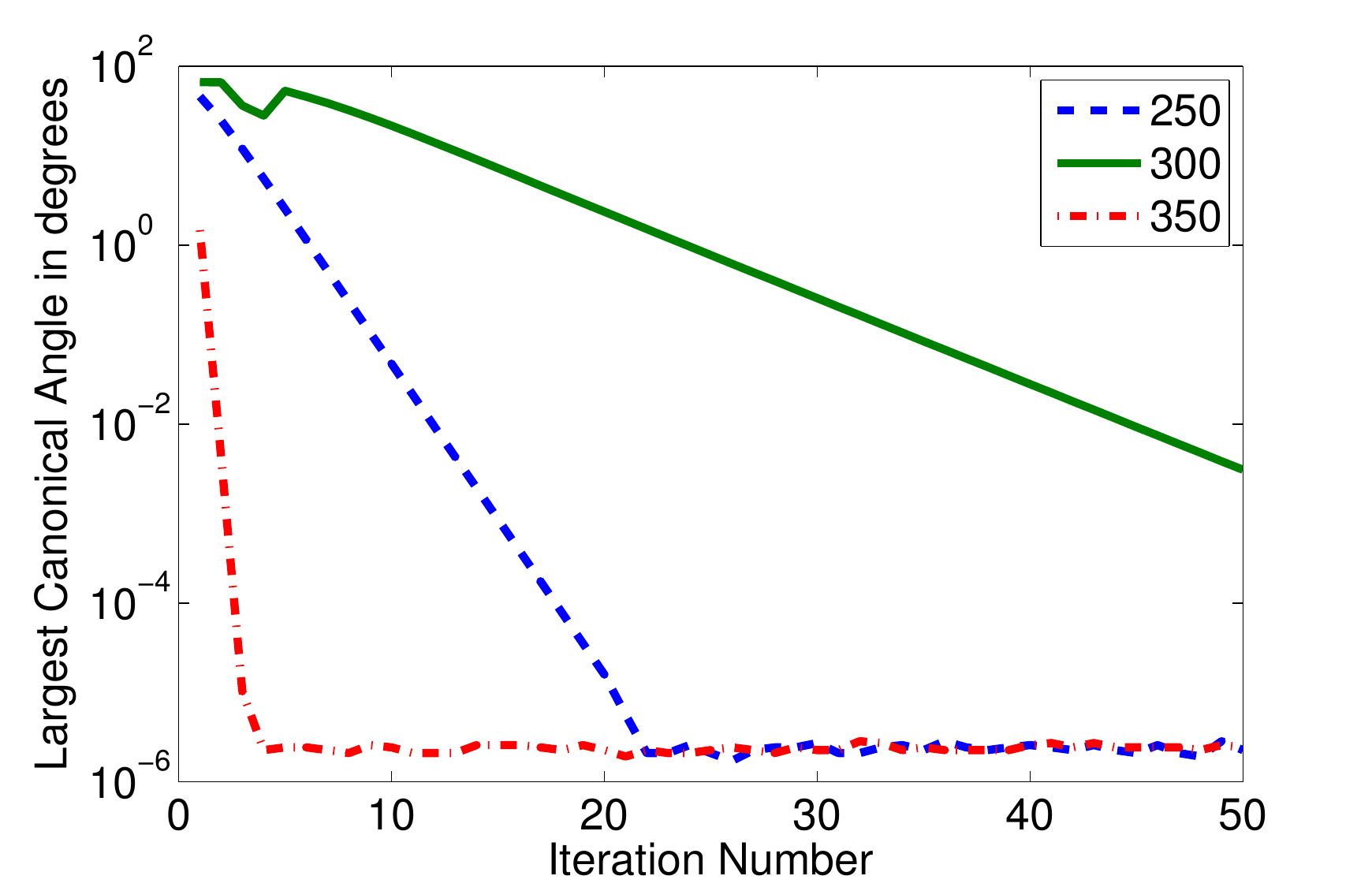}
      \quad
      \includegraphics[width=\figwidth]{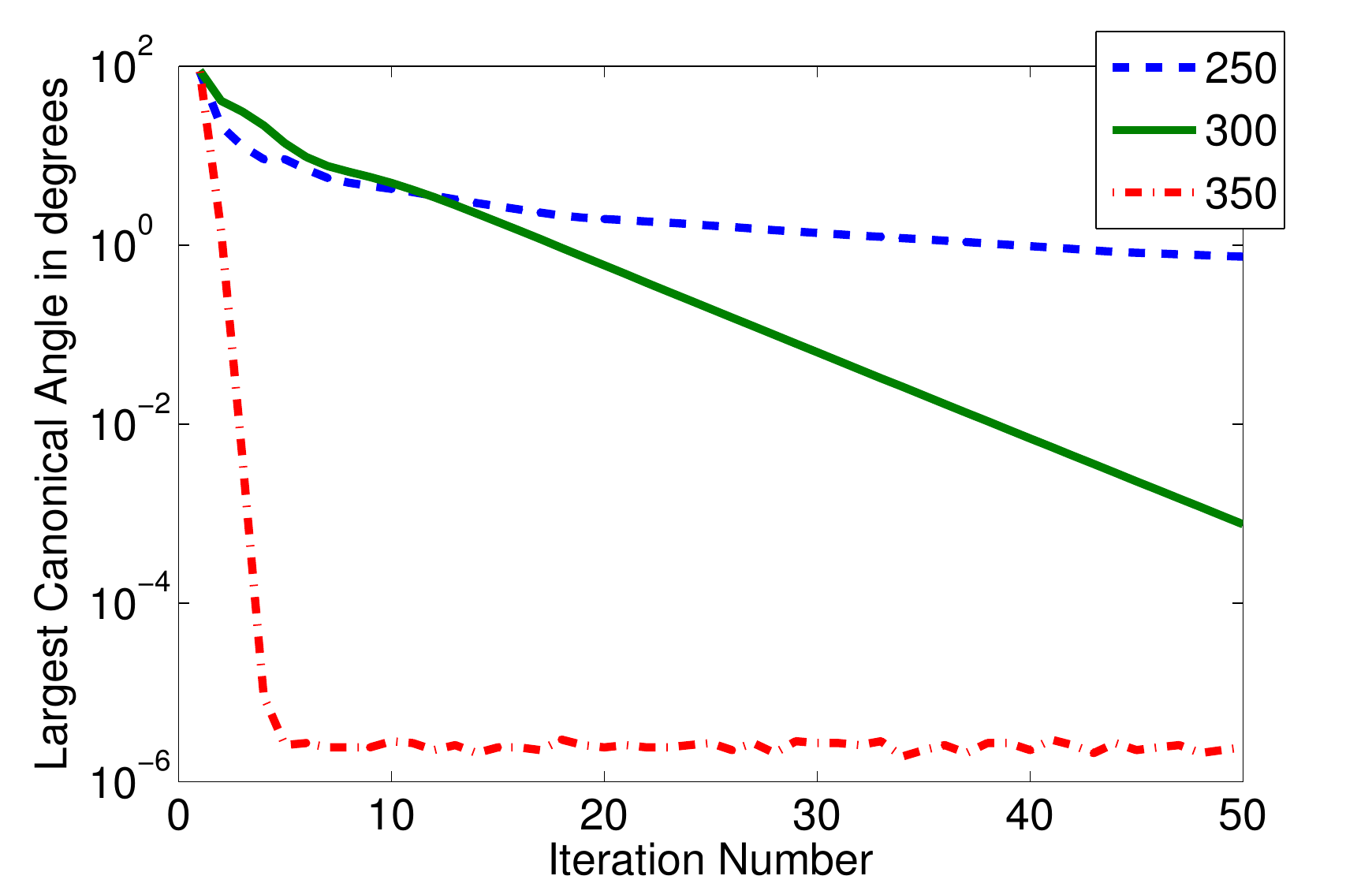}
   \end{center}
   \caption{
Canonical angles (left: between current iterate $X^{(i)}$ and exact eigenspace $X_I$; right: between current iterate $X^{(i)}$ and previous iterate $X^{(i-1)}$) in degrees for $\widetilde{M} = 250,300,350$.
           }
    \label{can_angle:fig}
\end{figure}

\begin{experiment}\label{searchspace:ex}
   We consider a size-$1059$ matrix $A = \texttt{LAP\_CIT\_1059}$ from
   modelling cross-citations in scientific publications, and $B = I$. In this test we search for
   $\widetilde{M}=1,\ldots,450$ eigenpairs with eigenvalues in an
   interval $\ilambda$ containing the $M=300$ lowest
   eigenvalues. The maximum number of iterations allowed for \feast is 20.

   The left panel of Figure~\ref{no_iters-and-minmaxres:fig} shows the number of iterations necessary for \feast to calculate all eigenpairs within $\ilambda$ with sufficiently small residual
   $\norm{Ax-\lambda Bx}\leq \varepsilon \cdot n \cdot \max\lbrace \abs{\lmin},\abs{\lmax}\rbrace$, as a function of $\widetilde{M}$.
   An iteration count of 20 typically implies that either none or not all eigenpairs converged within these 20 iterations.
   The right panel shows the residual span for all computed eigenpairs with eigenvalues in the interval after the respective number of iterations (20 or fewer, if convergence was reached beforehand). Again these numbers are given as a function of $\widetilde{M}$. We see that, leaving aside the very small region around the exact eigenspace size, either all or none of the eigenpairs show a sufficiently small residual.
   While for $\widetilde{M}<M$ no eigenpairs converge and especially the minimum residuals are large, for $\widetilde{M}>M$ also the maximum residuals begin to drop significantly and typically all eigenpairs may converge if only enough iterations are performed. With $\widetilde{M}$ just slightly larger than $M$, all eigenpairs reach convergence within few iterations.

For a better understanding of the evolution of the computed eigenspace, we monitored the largest canonical angle $\sphericalangle \left( X^{(i)},X_{\ilambda} \right)$ \cite[p.~603]{golubvanloan:1996} between the current approximate eigenspace $X^{(i)}$ and the exact eigenspace $X_{\ilambda}$, as well as the angle $\sphericalangle \left( X^{(i)},X^{(i-1)} \right)$ between the current and the previous iterate. Figure~\ref{can_angle:fig} provides these angles for three values of $\widetilde{M}$, $\widetilde{M}=250,\; \widetilde{M}=300$ and $\widetilde{M}=350$. In this last case, after five iterations the computed eigenspace contains the exact one and does not vary anymore; these two facts imply convergence. By contrast, the curves for $\widetilde{M}=250$ indicate that while the computed eigenspace becomes contained in the exact one after more than 20 iterations, it keeps varying, never to reach convergence. 
Interestingly, the worst convergence with respect to the exact eigenspace seems to occur for $\widetilde{M}=300$. This can be intuitively understood by the fact that two subspaces of the same dimension need to be identical in order to have an angle of zero between each other.
\end{experiment}

\subsection{Choice of the starting basis $Y$}\label{sec:StartingBase}

The choice of the starting basis
$Y \in \mathbb{C}^{n\times\widetilde{M}}$ in line~\ref{integral:line}
of Algorithm~\ref{basicfeast:alg} plays a critical role in the unfolding of the algorithm:
most importantly, $Y$ has to include components along $U$, so that its projection through $Q$ spans $\calx$.
Vice versa, if one or more of the columns of $Y$ are $B$-orthogonal to the space $\mathcal{U}$,
then the corresponding columns of $U=QY$ will be zero. 
If a good initial guess for the eigenvectors of $(A,B)$ is available, 
then it can used as the starting base 
$Y$; otherwise the typical choice is a base of random vectors.

\begin{experiment}
   We used \feast with $\widetilde{M}=450$ to compute the $M=300$ lowest
   eigenvalues and corresponding eigenvectors of the matrix pencil
   $( A = \texttt{LAP\_CIT\_1059}, B = I)$.
   These eigenvalues are simple and sufficiently far away from
   the only multiple eigenvalue, which is zero.
   With a fixed random starting basis $Y$, four iterations were
   sufficient to compute all wanted eigenpairs with residuals
   $\norm{Ax_j-\lambda_jBx_j} \leq 5.5\times 10^{-15}$.
   Then we projected out the ten eigenvectors corresponding to the
   $10$ lowest eigenvalues via
   $Y := ( I - X_{1 : 10} X_{1 : 10}^T) \cdot Y$.
   It took seven iterations for the lowest $290$ eigenpairs to
   converge, and nine more iterations for other nine.
   One eigenpair did not converge within the limit of 20 iterations. Convergence took place ``from top to bottom'', \ie the eigenpairs $11,\ldots,300$ converged first, then the eigenpairs $2,\ldots,10$. The smallest eigenvalue did not converge within the iteration limit.
\end{experiment}

In general, thanks to round-off errors, convergence could still
be reached in most of our tests. In fact even though 
some components were zeroed out, the floating point arithmetic causes 
almost zero entries to grow as the computation progresses. 
Convergence is then reached with noticeably more iterations.

\subsection{Stopping criteria}\label{stopping_crit:sec}

Algorithm 1 relies on a stopping criterion to determine whether the eigenpairs
are computed to a sufficient degree of accuracy. Such a criterion must balance
cost and effectiveness. In the original implementation \cite{polizzi:2009}, \feast
monitors convergence through the change in the 
\emph{sum} of the computed eigenvalues; more precisely, a relative criterion of the form 
\begin{equation}
	\label{trace_crit:eq}
   \frac{\abs{\mathsf{trace}_k-\mathsf{trace}_{k-1}}}{\abs{\mathsf{trace}_k}}
   <
   \mathrm{TOL}
\end{equation}
is used, where $\mathsf{trace}_k$ denotes the sum of the Ritz values
in the $k$th iteration lying in the search interval and $\mathrm{TOL}$ is a user tolerance.
Criterion \eqref{trace_crit:eq} raises three problems.

First, its denominator might be zero or close to, causing severe numerical instabilities. 
Such a scenario arises, for instance, when all the eigenvalues in the search interval are zero.
Second, the number $\mathsf{trace}_k$ can be (almost) constant for two
consecutive values of $k$, stopping \feast even if the residuals are still large.
The third problem is of a more general nature: 
if the algorithm stagnates before the considered eigenpairs converge,
all criteria that are based only on the change in the eigen\emph{values} 
might still signal convergence.
It is not hard to construct examples where this happens \cite{stewart:2001}.

As an alternative criterion we propose a per-eigenpair residual that depends on
the search interval: a Ritz pair has converged if it fulfills the inequality
\begin{equation}
	\label{eq:resy}
   \norm{Ax-Bx\lambda}
   \leq
   \varepsilon \cdot n \cdot \max\lbrace \abs{\lmin},\abs{\lmax}\rbrace
   .
\end{equation}
The extra cost of this criterion is \emph{one} matrix--vector product per
vector because $Bx$ is needed anyway to compute $Y$ in the next
iteration, or, if sparsity is not exploited,
$\mathcal{O}(\widetilde{M}^{2} n)$ operations since $Ax = (AU)w$,
with $A U$ reused from Step~2 of the \feast algorithm.
If $\max\lbrace \abs{\lmin},\abs{\lmax}\rbrace$ is too small, \eg when
only eigenpairs with zero eigenvalues are sought after, one may replace this
quantity by an estimate for $\norm{B^{-1}A}$, which is the magnitude of the largest eigenvalue of the problem and might be obtained by some auxiliary routine (\eg by some steps of a Lanczos method, see \cite{parlett:1998}).

Let us consider a matrix with spectrum symmetric with respect to zero, and choose the search interval to be symmetric around zero. If summed, the computed Ritz values cancel themselves pairwise, thus $\mathsf{trace}_k$ is approximately zero.
When \feast is applied to such a matrix, the trace criterion signals convergence when the difference in \eqref{trace_crit:eq} comes close enough to zero;
in this case, the ratio \eqref{trace_crit:eq} holds no information about convergence.
From multiple tests with random starting bases, we experienced that it took between 7 and 100 iterations for criterion \eqref{trace_crit:eq} to signal convergence. By contrast, criterion \eqref{eq:resy} dropped below $10^{-15}$ always after 6 iterations.

In an additional test, we ran \feast on a symmetric matrix of size $470$, seeking the known $57$-fold eigenvalue $1$. Here we chose the search interval symmetrically around $1$, and $\widetilde{M}=120 > M$.
After  five iterations, \eqref{trace_crit:eq} was $1.2\times 10^{-16}$, effectively halting the computation, although the residuals \eqref{eq:resy} were still of order $10^{-9}$ to $10^{-12}$. The right hand side of \eqref{eq:resy} was about $10^{-13}$ in this example, meaning that none of the residuals was satisfactory small. 

\subsection{The impact of the linear solver
            on residuals and orthogonality}\label{lin_solvers:sec}

As seen in Section~\ref{sec:IntegratingTheResolvent}, the computation
of the basis according to \eqref{eq:u-integral} involves solving
several linear systems of the form
\begin{equation}
 (zB-A)V = B Y
 \label{linsys:eq} 
\end{equation}
for $V$.
The particular values of the integration points $z$ depend on the
method chosen for numerical integration, which is not discussed here.
Typically, $z$ will be a complex number near the spectrum of $(A,B)$.
Recall that \eqref{linsys:eq} is a linear system with $\widetilde{M}$ right
hand sides. In principle, any linear solver can be used.
Direct solvers, \eg Gaussian elimination based, can be prohibitively
expensive because for each value of $z$ we need to factorize $zB-A$,
which may be an $\mathcal{O}(n^3)$ process. 

The methods of choice for solving large sparse linear systems without
further knowledge about the underlying problem are Krylov subspace
methods; for a review see, \eg \cite{saad:2003}.
Here we cannot give a detailed discussion, but let us remark that the
convergence of Krylov subspace methods depends on several parameters.
First, the best convergence results can be expected for Hermitian
matrices since this property can be exploited.
Unfortunately, \eqref{linsys:eq} typically has a non-real diagonal and
therefore is not a Hermitian problem.
(However, if $B = I$ then the matrix $z I - A$ is \emph{shifted}
Hermitian, so methods for shifted systems may be applicable
\cite{frommer_glaessner:1998}.)
Second, the convergence for a fixed method typically depends on the
structure of the spectrum of the matrix.
The eigenvalues of $zB-A$ are scattered over the complex plane so that
no good convergence results can be inferred.
Third, the condition number $\norm{(zB-A)^{-1}}\cdot\norm{zB-A}$ of the
system plays an important role and is often large for
\eqref{linsys:eq}, since $z$ can be very close to the spectrum of
$(A,B)$.
For these reasons, standard Krylov subspace solvers may need a large
number of iterations to converge.
This expectation was confirmed by our experiments.
The need for an effective preconditioner is apparent, and its
development is part of further research.

Another way to speed up the linear solvers is to terminate them before
full convergence is reached.
Thus the question arises how
accurately the systems \eqref{linsys:eq} need to be solved in order to
obtain eigenpairs of sufficient quality
in a reasonable number of \feast iterations.
We therefore investigated the effect of the accuracy in the
solution of linear systems on the ultimately achievable
per-eigenpair residuals and the orthogonality of the eigenvectors, as
well as on the number of \feast iterations.

\begin{figure}[htb]
  \begin{center}
    \includegraphics[scale=0.5]{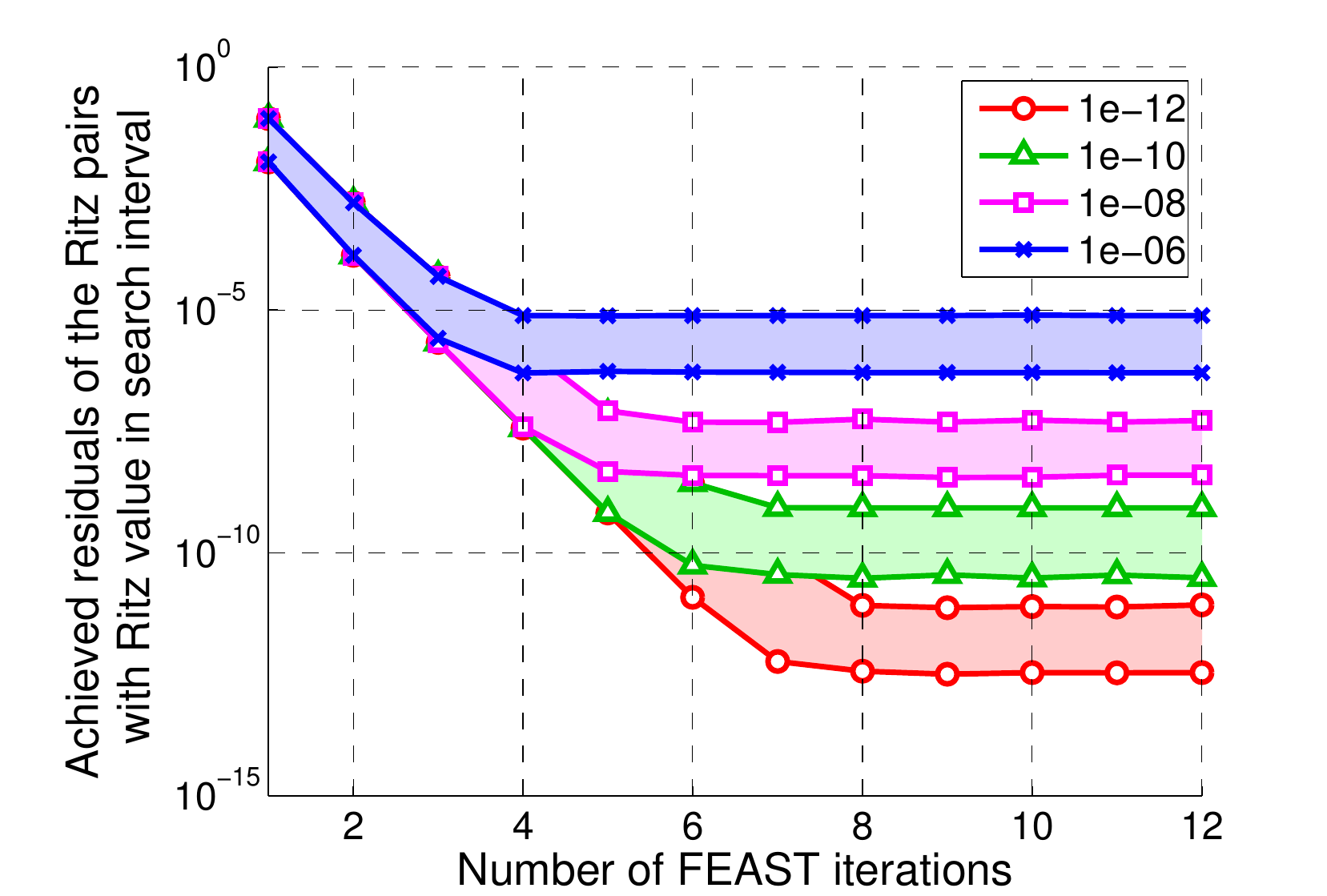}
  \end{center}
  \caption{Range of all residuals among all Ritz pairs 
           in $\ilambda$ for four different residual bounds
           $\varepsilon_{\mathrm{lin}}$ in the linear solver.}
  \label{linsolv:fig}
\end{figure}

\begin{experiment}
We applied Algorithm~\ref{basicfeast:alg} to the matrix pair
($A$,$B$), where $A = \texttt{LAP\_CIT\_395}$ arises
in the modelling of cross-citations in scientific publications, and $B$ was chosen to be a diagonal matrix with random entries.
We calculated the eigenpairs corresponding to the
$10$ largest eigenvalues.
The linear systems were solved column-by-column by
running GMRES \cite{saad:2003} until
$\norm{(zB-A)v_j -By_j}/\norm{r^0_{j}} \leq \varepsilon_{\mathrm{lin}}$,
where $r^0_j$ is the starting residual.
Figure~\ref{linsolv:fig} reveals that the residual bounds that were
required in the solution of the inner linear systems translated almost
one-to-one into the residuals of the Ritz pairs.
Even for a rather large bound such as
$\varepsilon_{\mathrm{lin}} = 10^{-6}$, the \feast algorithm still
converged (even though to a quite large residual).
For the orthogonality of the computed eigenvectors $x_{j}$, the
situation was different.
After $20$ \feast iterations, an orthogonality level
$\max_{i\neq j}\abs{x_i^HBx_j}$ of order $10^{-15}$ could be reached for
each of the bounds
$\varepsilon_{\mathrm{lin}} = 10^{-6}, 10^{-8}, 10^{-10}, 10^{-12}$
in the solution of the linear systems.
Thus the achievable orthogonality does not seem to be
very sensitive to the accuracy of the linear solves.
It also did not deteriorate significantly for a
larger number of desired eigenpairs.
\end{experiment}

\section{\feast for multiple intervals }\label{feast_in_context:sec}

Using \feast for computing a large number $M$ of eigenpairs is not
recommended since the complexity grows at least as
$\mathcal{O}( n M^{2} )$
due to the matrix--matrix products in Steps~2 and 4 of
Algorithm~\ref{basicfeast:alg}, or $\mathcal{O}( n^{2} M )$
if sparsity is not exploited, because $\mathcal{O}( M )$
matrix--vector products must be computed.
(In fact, if $M$ approaches $n$ then just Step~3 of
Algorithm~\ref{basicfeast:alg} has roughly the same complexity as
the computation of the full eigensystem of the original problem.)
However, as already hinted at in \cite{polizzi:2009},
\feast's ability to determine the eigenpairs in a specified
interval makes it an attractive building block to
compute the eigenpairs by subdividing the search interval $\ilambda$ into
$K$ subintervals $I_{\lambda}^{(1)}, \ldots, I_{\lambda}^{(K)}$ 
and applying the algorithm---possibly in parallel---to each one of them.

In the following we report on issues concerning the orthogonality of
eigenvectors coming from different subintervals.
It is well known that eigenvectors computed
independently from each other tend to have worse orthogonality than
those obtained in a block-wise manner.
Furthermore, the quality of the results depends on the internal structure of the
spectrum, namely the relative distances between the eigenvalues.
For more details, see, \eg \cite{parlett:1998}.

In the following we distinguish between {\em global} and {\em local} orthogonality,
according to the definitions  
\[
   \mathsf{orth}_{\mathrm{global}}
   = \max_{i\neq j, \; \lambda_i, \lambda_j \in I_{\lambda}} \abs{x_i^HBx_j}
   \quad
   {\rm and}
   \quad
   \mathsf{orth}_{k}
   = \max_{i\neq j, \; \lambda_i, \lambda_j \in I_{\lambda}^{(k)}} \abs{x_i^HBx_j}.
\]
Note that $\mathsf{orth}_{\mathrm{global}}$
denotes the worst orthogonality among \emph{all} computed eigenvectors,
while $   \mathsf{orth}_{k}$
describes the orthogonality achieved locally for the $k$th subinterval $I_{\lambda}^{(k)}$.   
The next two experiments reveal quite different behavior, depending
on the presence of clusters and on the choice of subintervals.

\begin{figure}[htb]
   \begin{center}
      \includegraphics[scale=0.5]{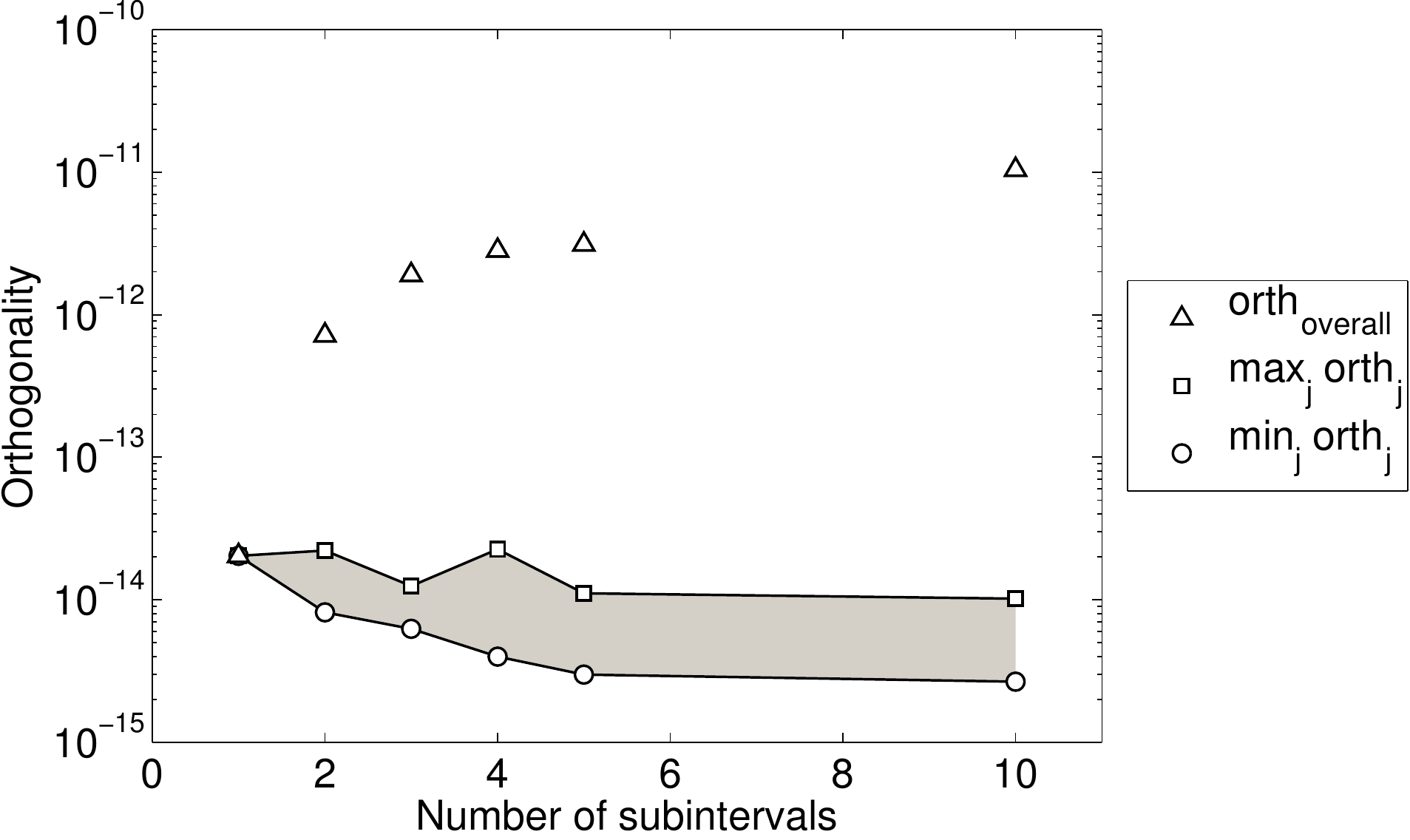}
   \end{center}
   \caption{Global orthogonality and range of local orthogonality $\mathsf{orth}_{j}$ for $K = 1,\ldots,5$ and $K=10$.}
   \label{orth:fig}
\end{figure}

\begin{experiment}\label{1473:ex}
In this experiment we calculate the $800$ lowest eigenpairs of the size-$1473$ matrix
pair (\texttt{bcsstk11}, \texttt{bcsstm11}) from the Matrix Market
(\texttt{http://math.nist.gov/MatrixMarket/}).
The corresponding eigenvalues range from $10.5$ to $3.8 \times 10^{7}$ and are not
clustered. We utilize different numbers of subintervals, $K=1,\ldots,5$, and
$K = 10$.
Figure~\ref{orth:fig} shows that while the local orthogonality is high and maintained
as the number of intervals increases, the global one 
degrades by two or more orders of magnitude. 
\end{experiment}

\begin{figure}

   \begin{center}
	\begin{tabular}{ccc}
         \includegraphics[width=60mm]{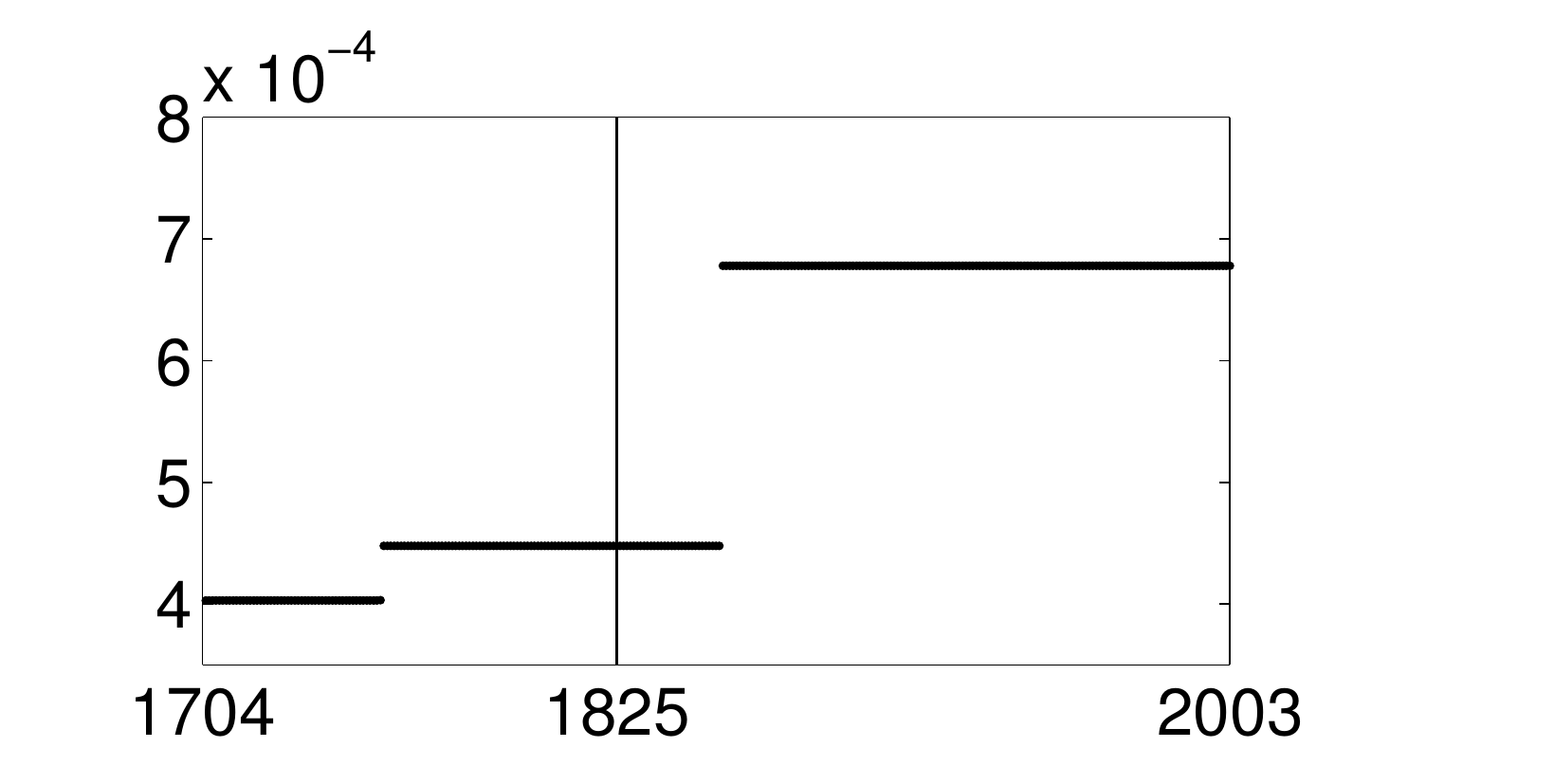}
	   & \qquad &
         \includegraphics[width=60mm]{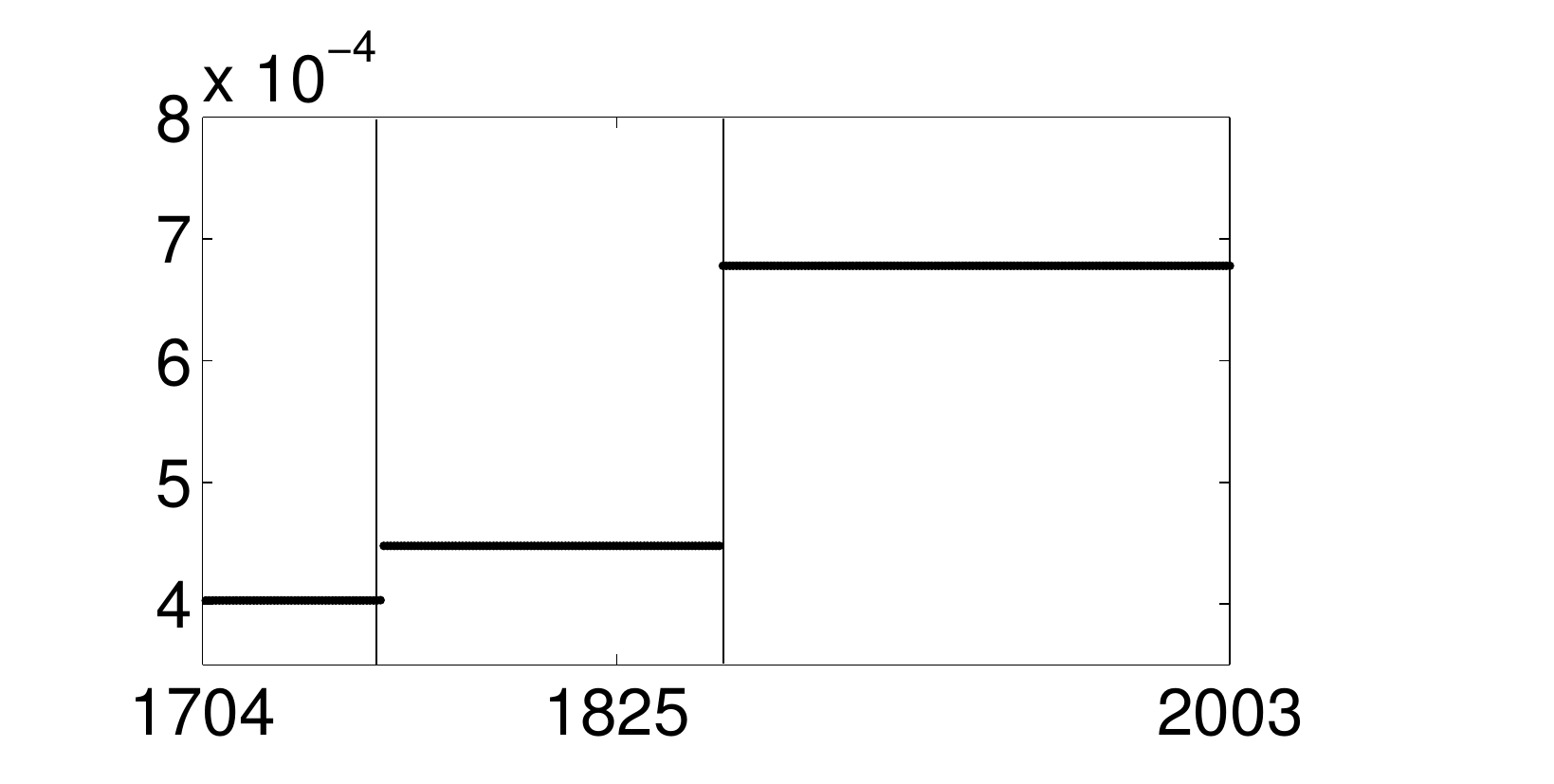}
	   \\
         \includegraphics[width=60mm]{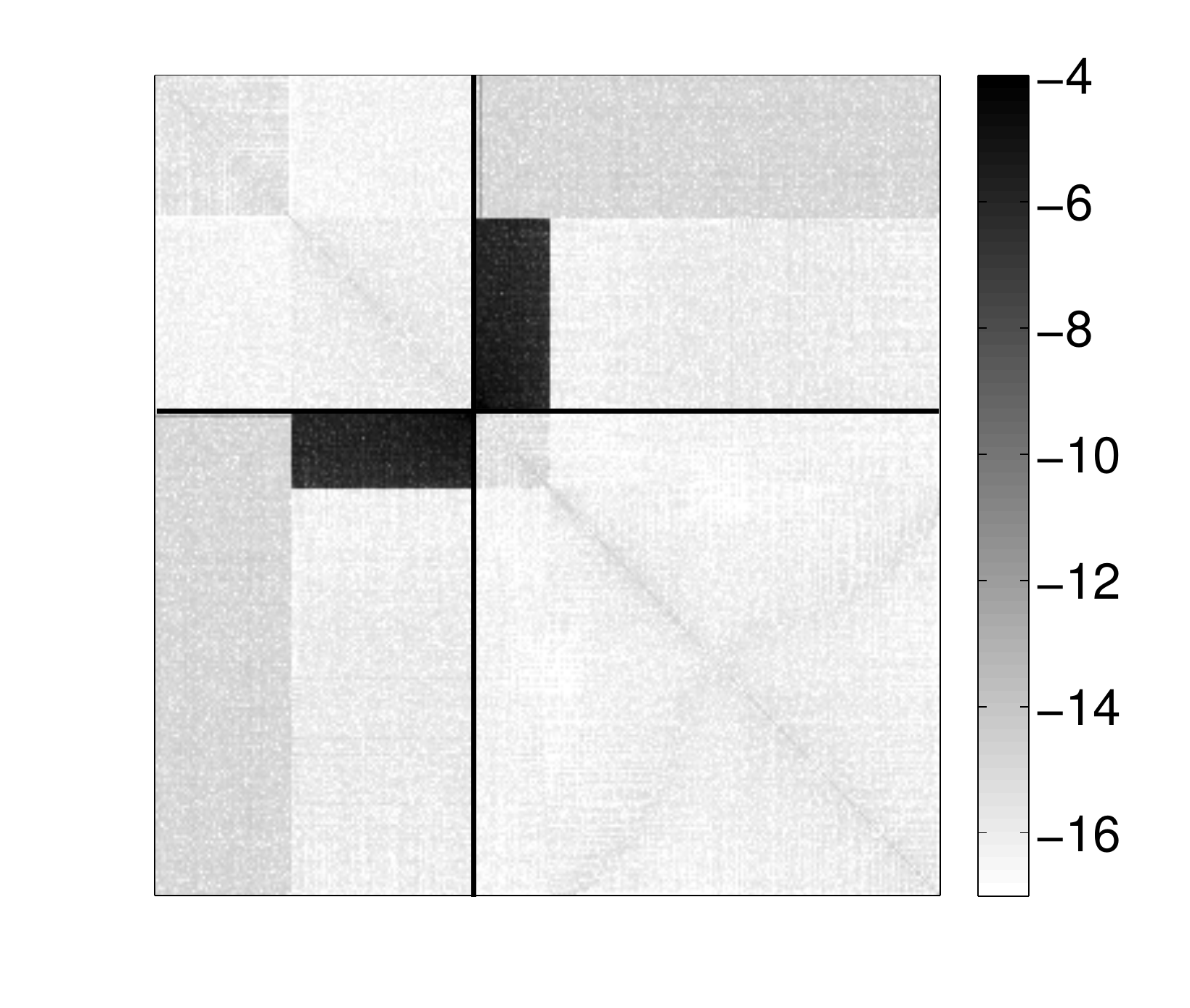}
	   & \qquad &
         \includegraphics[width=60mm]{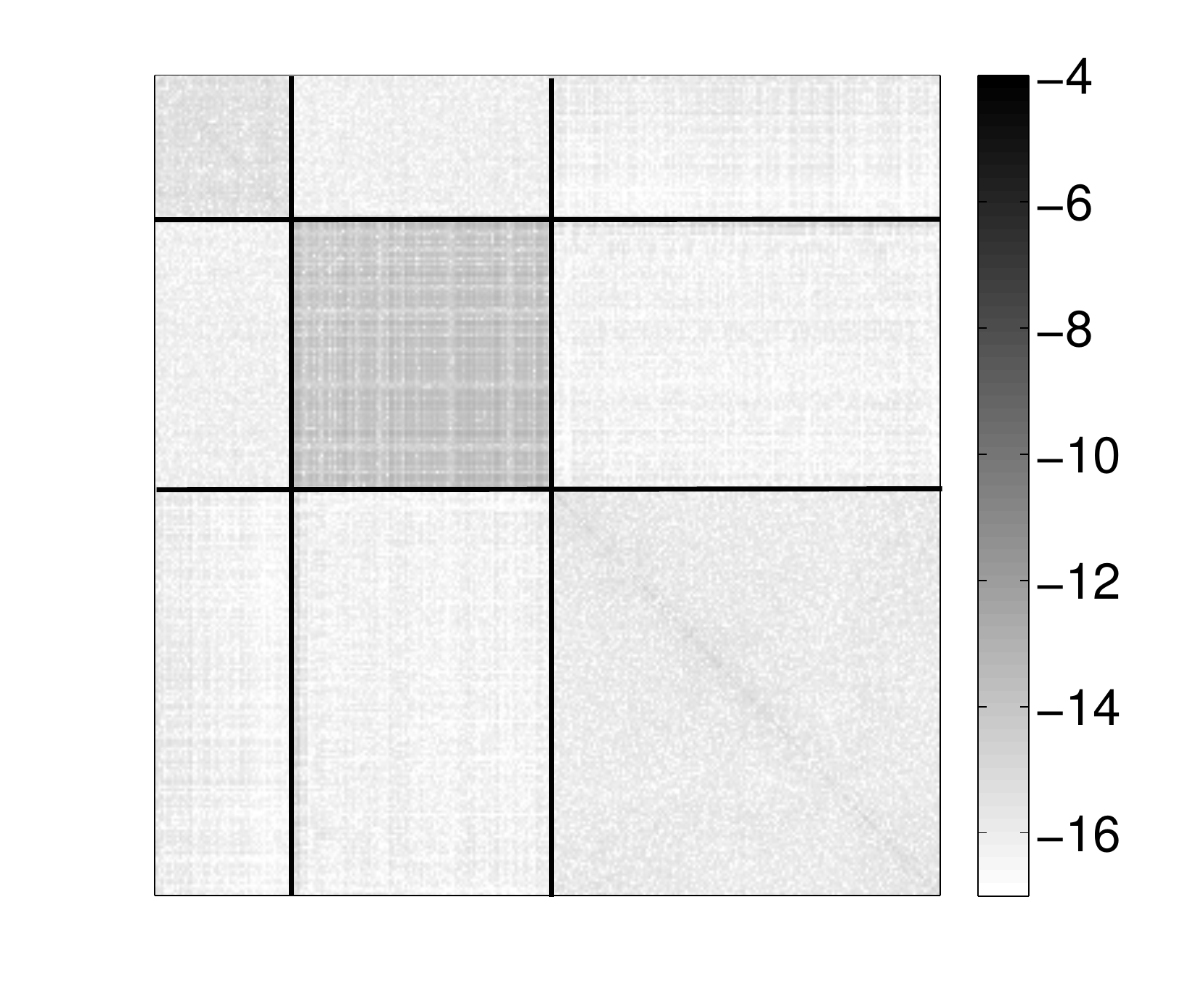}
	   \\[-1em]
	\end{tabular}
   \end{center}

    \caption{Computation of the eigenpairs corresponding to the $300$
             largest eigenvalues $\lambda_{1704}, \ldots, \lambda_{2003}$
             with the subdivision point $\mu = \lambda_{1825}$ taken from
             a group of very close eigenvalues.
             The top plots show the eigenvalues and the subdivision points
             (vertical lines), the bottom plots give a pictorial visualization of the orthogonality
             $\abs{x_{i}^{H} B x_{l}}$, $i \not= l$.}
    \label{orth_level:fig}
\end{figure}

\begin{experiment}\label{2003:ex}
In this test we consider a real unreduced tridiagonal matrix $A$ of size $2003$.
Its eigenvalues are simple, even though some
are tightly clustered; see the top plots in
Figure~\ref{orth_level:fig}.
The objective is to compute the 300 largest eigenpairs. 
To this end we initially split the interval
$\ilambda = \left[\lambda_{1704},\;\lambda_{2003}\right]$
into $\ilambda^{(1)} = [ \lambda_{1704}, \mu ]$
and $\ilambda^{(2)} = [ \mu, \lambda_{2003} ]$, with
$\mu = \lambda_{1825} \approx 0.448 \times 10^{-3}$ chosen within a
a cluster of $99$ eigenvalues. The relative gap between eigenvalue $\lambda_{1825}$
and its neighbors is of about $10^{-12}$
 (\ie agreement to roughly eleven leading decimal digits).
A sketch of the eigenspectrum with $\mu$ is given in the top left
of Figure~\ref{orth_level:fig}.
  
While \feast attains very good local orthogonality for both subintervals
 ($\mathsf{orth}_{1} = 4.4\times 10^{-15}$ and 
$\mathsf{orth}_{2} = 5.7\times 10^{-14}$), it fails to deliver global orthogonality
($4.7\times 10^{-4}$).
In the bottom left plot of Figure~\ref{orth_level:fig}
we provide a pictorial description of $\abs{x_{i}^{H} B x_{j}}$, $\lambda_i, \lambda_j \in \ilambda$.
The dark colored regions indicate that the loss of orthogonality 
emerges exclusively from eigenvectors belonging to the 99-fold cluster.
Next we divide the interval into 3 segments making sure not to break
existing clusters (see top right of Figure~\ref{orth_level:fig}). 
As illustrated in the 
bottom right plot, both the local and global orthogonality are satisfactory ($10^{-13}$ or better).
\end{experiment}

\section{Conclusions}\label{sec:concl}

We expounded the close connection between \feast algorithm and the
well-established \RRm for computing selected eigenpairs of a
generalized eigenproblem.  Starting from the mathematical foundation of this connection,
we identified aspects of the solver that might play a critical role
for its accuracy and reliability. Specifically, we discussed the 
choice of starting basis and stopping criterion, and 
the relation between the accuracy of the solutions of the linear systems internal to \feast 
and the resulting eigenpairs; we also investigated the use of \feast for computing a large portion or even the entire eigenspectrum.  
Through numerical examples we illustrated how each of these aspects might affect the robustness of the algorithm or
diminish the quality of the computed eigensystem.

While we hinted at possible improvements for several of the existing
issues, some questions remain open and are the subject of further
research.  For instance, a mechanism is needed to overcome the problem
of having to specify both the boundaries of the search interval and the
number of eigenvalues expected. Additionally, it would be desirable to
have a flag assessing whether all the eigenpairs in the search
interval have been found. Finally, orthogonality across multiple
intervals should be guaranteed.
 
In summary, our findings suggest that at the moment \feast is a
promising eigensolver for a certain class of problems, \ie when
a small portion of the spectrum is sought and knowledge of the
eigenvalue distribution is available. On the other hand, we believe it
is still not yet competitive as a robust ``black~box'', general-purpose
solver\footnote{At the time of submission, v2.0 of the \feast Package Solver was announced. This new version is designed for parallel computation, and doesn't seem to address the issues raised in this research paper. We still would like to remark that all the numerical experiments in this paper were executed making use of the algorithm from \feast v1.0.}.

\section*{Acknowledgements} The authors want to thank Mario Th\"une from MPI MIS in Leipzig for providing the \texttt{LAP\_CIT} test matrices.


\bibliographystyle{elsarticle-num}
\bibliography{feast}

\end{document}